\documentclass[a4paper,10pt]{article}

\usepackage{amssymb,amsmath,graphics}
\usepackage{enumerate}
\usepackage{bbm}
\newenvironment{dwd}{\par\noindent{\bf Proof.}}{\par\rightline{$\blacksquare$}}

\newtheorem{theo}{Theorem}
\newtheorem{prop}{Proposition}  
\newtheorem{coro}{Corollary}

\newtheorem{conj}{Conjecture}
\newtheorem{rema}{Remark}
\newtheorem{lema}{Lemma}
\newtheorem{defi}{Definition}
\def\be#1\ee{\begin{equation}#1\end{equation}}
\newcommand{\ba}{\begin{eqnarray} }
\newcommand{\ea}{\end{eqnarray} }
\def\bt#1\et{\begin{theo}#1\end{theo}}
\def\bl#1\el{\begin{lema}#1\end{lema}}
\def\bp#1\ep{\begin{prop}#1\end{prop}}
\def\bd#1\ed{\begin{defi}#1\end{defi}}
\def\ccA{{\cal A}}
\def\ccB{{\cal B}}

\def\va{\varepsilon}
\def\ra{\rightarrow}

\def\E{\mathbf{E}}

\def\P{\mathbf{P}}

\def\N{{\mathbb N}}

\def\R{{\mathbb R}}
\def\Z{{\mathbb Z}}

\def\ls{\leqslant}
\def\gs{\geqslant}

\begin{document}

\title{\bf On a contraction property of Bernoulli canonical processes}
\author{{Witold Bednorz \& Rafa\l{} Martynek  }
\footnote{Research partially supported by  MNiSW Grant N N201 608740.}
\footnote{Institute of Mathematics, University of Warsaw, Banacha 2, 02-097 Warszawa, Poland}}
\date{}

\maketitle

\begin{abstract}
In this paper we give several results concerning the supremum of canonical processes. The main theorem concerns a contraction property of Bernoulli canonical process which generalizes the one proved by Talagrand (Theorem 2.1 in \cite{Tal2}). The result works for independent Rademacher random variables $(\va_i)_{i\geq1}$ and states that we can compare $\E\sup_{t\in T}\sum_{i\geq1}\varphi_{i}(t)\va_i$ with  $\E\sup_{t\in T}\sum_{i\geq1}t_i\va_i$, where a function $\varphi=(\varphi_i)_{i\geq1}: \ell^2\supset T\rightarrow\ell^2$, satisfies certain conditions. Originally, it is assumed that each of $\varphi_i$ is a contraction. We relax this assumption towards comparison of Gaussian parts of increments, which can be described in the following way. For all $s,t\in T$, $p\gs 0$
$$
\inf_{|I^c|\ls Cp}\sum_{i\in I}|\varphi_i(t)-\varphi_i(s)|^2\ls C^2\inf_{|I^c|\ls p}\sum_{i\in I}|t_i-s_i|^2,
$$
where $C\gs 1$ is an absolute constant and $I\subset\N$, $I^c=\N\backslash I$.
\end{abstract}

\section{Introduction and notation}

Throughout this paper we will use the following notation. For the set $A$ the number of elements in $A$ will be denoted as $|A|$. If $t=(t_n)$, ${n\gs1}$ is a sequence of real numbers and $p\gs1$ then $\|t\|_p=(\sum_{n=0}^{\infty}t_i^p)^{\frac{1}{p}}$ and $\ell^p$ is the space of all sequences $t$ with $\|t\|_p<\infty$. If $S,T\subset\ell^p$ then $S+T=\{s+t:s\in S, t\in T\}$. For a random variable $\xi$ and $p>0$ we put $\|\xi\|_p=(\E|\xi|^p)^{\frac{1}{p}}$.
If $(\xi_i)$, $i\gs 1$ is a sequence of independent, identically distributed random variables such that $\E \xi_i=0$, $\E \xi_i^2=1$
and $t=(t_n)\in\ell^2$ then the random variable

\be\label{repres}
X_t=\sum^{\infty}_{i=1}t_i \xi_i,
\ee
is well-defined. For each $T\subset\ell^2$ with $0\in T$ the process $X_T=(X_t)_{t\in T}$ is called \textit{canonical}. The convergence of the above series holds in the sense of $\|\cdot\|_2$ which means that
$$
\lim_{n\ra\infty}\|\sum^n_{i=1}t_i\xi_i- X_t\|_2=0.
$$ 
Clearly,
$$
\|X_t-X_s\|_2=\|t-s\|_2,\;\;\mbox{for}\;\;s,t\in T.$$
\begin{rema}
The almost sure convergence in (\ref{repres}) might be guaranteed also when the independence assumption on $\xi_i$'s is skipped. In such case we may consider finite dimensional version of (\ref{repres}), where $T\subset \R^d$. The most studied example is when $\xi_i$'s have log-concave tails i.e. $\P(|\xi_i|>t)=\exp(-N_i(t))$ for $N_i:[0,\infty]\rightarrow[0,\infty]$ convex and may be dependent.
\end{rema}
We want to distinguish two types of canonical processes which will be of special interest. If $(\xi_n)=(\va_n)$ and $\P(\va_n=1)=\P(\va_n=-1)=\frac{1}{2}$ then the process $X_T$ is called \textit{canonical Bernoulli} and it is denoted by $B_T=(B_t)_{t\in T}$. This class of processes is important for various applications
e.g. infinitely divisible processes \cite{Tal2}, empirical processes (see \cite{Tal1} for the comprehensive study). If $(\xi_n)=(g_n)$ and $g_n$ are distributed by the normal law $\mathcal{N}(0,1)$ then the process $X_T$ is called \textit{canonical Gaussian} and it is denoted by $G_T=(G_t)_{t\in T}$. In fact, canonical Gaussian processes can be seen as a motivation to study canonical processes in general. The reason for that being the Karhunen-L\`{o}eve representation of separable Gaussian process with the canonical Gaussian process, (see e.g. \cite{M-R} Corollary 5.3.4).\\
\noindent The main object studied will be the suprema of canonical processes. For any set $T$ and a stochastic process $(X_t)_{t\in T}$ we define 
$$
S_X(T)=\sup_{F\subset T}\E\sup_{t\in F}X_t,
$$
where the supremum is taken over all finite subsets $F$ of $T$.
Usually, by considering the separable modification of $X_t, t\in T$ it is possible to guarantee that $\sup_{t\in T}X_t$ is well-defined random variable (for the definition of separable version of the process and the discussion on the measurability of the supremum in a general setting of Banach space which is not necessarily separable see Ch. 2 in \cite{L-T}).  In this case $S_X(T)$ coincides with the usual expectation of the supremum over $X_t$, namely
$$
S_X(T)=\sup_{F\subset T}\E\sup_{t\in F}X_t=\E\sup_{t\in T}X_t.
$$
Let us finish this section with a few important technicalities which will be helpful in dealing with canonical processes. We have that $S_X(T)=S_X(T-t)$, where $T-t=\{s-t:\;s\in T\}$
so we may always require that $0\in T$. Moreover, $S_X(T)=S_X(\mathrm{Conv}T)$
and $S_X(T)=S_X(\mathrm{cl}T)$, where $\mathrm{Conv}T$ is the convex 
hull of the set $T$ and $\mathrm{cl}T$ is the closure of $T$ in $\ell^2$.\\
We follow the convention that numerical constants denoted by the same letter might vary from line to line. The same constants will be subindexed e.g. $C_1, C_2$ etc.

\section{Suprema of canonical processes via chaining}

\noindent
First, we recall the basics of the chaining approach to upper bounds for stochastic processes. We say that the sequence $\ccA=(\ccA_n)_{n\gs 0}$ of partitions of $T$ is an \textit{admissible partition of T} if $\ccA_{0}=\{T\}$ and $|\ccA_n|\ls N_n=2^{2^n}$ for $n\gs 1$
(usually it is required also that these partitions are \textit{nested} i.e. for any set $A\in A_n$, $n\gs1$ there is a set $B\in A_{n-1}$ such that $A\subset B$). For $t\in T$ we denote by $A_n(t)$ the unique element of partition $\ccA_n$ which contains $t$.
Let $\pi_n:T\rightarrow T$ be a sequence of points such that its' $n$-th element is defined so that $\pi_n(t)=\pi_n(s)$ for all $s,t\in A\in \ccA_n$. We will denote it by $\pi_n(A)$. Let
$$
\gamma_X(T)=\inf\sup_{t\in T}\sum_{n\gs 1}\|X_{\pi_{n}(t)}-X_{\pi_{n-1}(t)}\|_{2^n},
$$
where the infimum is taken over all admissible sequences of partitions. We denote by $T_n$ 
the family of all $\pi_n(t)$, $t\in T$. In words, at each step of partitioning we choose some point $\pi_n(t)$ which belongs to the same partition set as t. Clearly, $|T_n|\ls 2^{2^n}$ and $|T_0|=1$. By the \textit{chain} we mean writing $X_t$ as the sum of consecutive approximations i.e.
$$X_t=X_{\pi_0(t)}+\sum_{n\gs1}(X_{\pi_{n}(t)}-X_{\pi_{n-1}(t)}).$$
Let us observe the following property of $\gamma_X(T)$.
\begin{lema}
Let $T_1$ and $T_2$ be some index sets. Suppose that for some stochastic process $(X_t)$ $\gamma_X(T_1)$, $\gamma_X(T_2)$ and $\gamma_X(T_1+T_2)$ are well-defined and for  $n\gs1$
\be\label{reg}
 \|X_t\|_{2^{n+1}}\ls C_1\|X_t\|_{2^n}.
\ee
Then,
\be\label{subadd}
\gamma_X(T_1+T_2)\ls C_1(\gamma_X(T_1)+\gamma_X(T_2)).
\ee
In particular, for canonical Bernoulli and canonical Gaussian processes the above inequality holds with $C_1=\sqrt{3}$.
\end{lema}
\begin{dwd}
Let $t\in T_1$ and $s\in T_2$. Let $(\ccA_n(T_1))_{n\gs0}$ and $(\ccA_n(T_2))_{n\gs0}$ be admissible partitions of $T_1$ and $T_2$ respectively. 
Define $\ccA_{n+1}(T_1+T_2)$ as all possible sums of these partitions i.e. $\ccA_{n+1}(T_1+T_2)=\{A+B: A\in\ccA_n(T_1), B\in\ccA_n(T_2)\}$ and $\ccA_0(T_1+T_2)=\{T_1+T_2\}$. It is obviously admissible since $N_n\cdot N_n=N_{n+1}$. Moreover, for $A\in \ccA_n(T_1)$ and $B\in \ccA_n(T_2)$ let $\pi_{n+1}(A+B)=\pi_n(A)+\pi_n(B)$. We also put $\pi_0(T_1+T_2)=\pi_0(T_1)+\pi_0(T_2)$. Clearly, for $t\in A$ and $s\in B$
$$
\|X_{\pi_{n+1}(t+s)}-X_{\pi_{n}(t+s)}\|_{2^{n+1}}\ls \|X_{\pi_{n}(t)}-X_{\pi_{n-1}(t)}\|_{2^{n+1}}+
 \|X_{\pi_{n}(s)}-X_{\pi_{n-1}(s)}\|_{2^{n+1}}.
$$  
So, by (\ref{reg}) we obtain
$$
 \|X_{\pi_{n+1}(t+s)}-X_{\pi_{n}(t+s)}\|_{2^{n+1}}\ls C_1(\|X_{\pi_{n}(t)}-X_{\pi_{n-1}(t)}\|_{2^{n}}+
 \|X_{\pi_{n}(s)}-X_{\pi_{n-1}(s)}\|_{2^{n}}).
$$
The conclusion follows since $t\in T_1$, $s\in T_2$ were arbitrary as well as the partitions $\ccA_n(T_1)$ and $\ccA_n(T_2)$. The reason why the inequality (\ref{reg}) holds with constant $\sqrt{3}$ for canonical Gaussian is a straightforward consequence of the fact that for $p\in\N$ we have $\|G_t\|_{2p}=\|t\|_2(2p)!!$, where $(2p)!!=(2p-1)(2p-3)\dots 1$. The result for canonical Bernoulli processes follows from the general Kahane's inequality (see e.g.\cite{Gar} Theorem 13.2.1)
\end{dwd}
In \cite{Lat},\cite{Me-P} it was proved that under a suitable regularity assumptions 
$S_X(T)\ls K\gamma_X(T)$, where $K$ is a universal constant. Let us give a short argument for a similar upper bound. 

\begin{theo}\label{theo2}
For a stochastic process $(X_t)_{t\in T)}$ for which $S_X(T)$ is well-defined we have that 
$$S_X(T)\ls 4\gamma_X(T)$$.
\end{theo}
\begin{dwd}
Let $(\ccA_n)_{n\gs0}$ be any admissible partition of $T$. For any set $A\in \ccA_n$ and $k\ls n$ we denote by $A_k$ its $k$-parent i.e. $A_k\in\ccA_k$
and $A\subset \ccA_k$. Consequently, if $t\in A\in \ccA_n$ then $\pi_k(t)=\pi_k(A_k)$. 
The proof is based on the analysis of the partition sequence.
Let $N$ be fixed and consider $\ccA_N$. The chaining argument gives
\begin{align*}
& \E\sup_{A\in \ccA_N}|X_{\pi_N(A)}-X_{\pi_0(A)}|\ls \E \sup_{A\in \ccA_N}\sum^N_{k=1}|X_{\pi_k(A_k)}-X_{\pi_{k-1}(A_{k-1})}|\\
&\ls \E \sup_{A\in \ccA_N}\sum^N_{k=1} 2\|X_{\pi(A_k)}-X_{\pi_{k-1}(A_{k-1})}\|_{2^k}(1+
\E \sup_{A\in \ccA_N}( \frac{|X_{\pi_k(A_k)}-X_{\pi_{k-1}(A_{k-1})}|}{2\|X_{\pi(A_k)}-X_{\pi_{k-1}(A_{k-1})}\|_{2^k}}-1)_{+}\\
& \ls 2\gamma_X(T)(1+ \sum^N_{k=1} \sum_{B\in \ccA_k}\E (\frac{|X_{\pi_k(B)}-X_{\pi_{k-1}(B_{k-1})}|}
{2\|X_{\pi_k(B)}-X_{\pi_{k-1}(B_{k-1})}\|_{2^k}}-1)_{+}).
\end{align*}
 We show that 
 $$
 \E(\frac{|X_t|}{2\|X_t\|_{2^k}}-1)_{+}\ls \frac{1}{2^kN_k}.
 $$
  Indeed, denoting $Y_t=X_t/\|X_t\|_{2^k}$ we have $\E |Y_t|^{2^k}=1$ and hence
\begin{align*}
\E((|Y_t|/2)-1)_{+}&=\int^{\infty}_{0}\P(|Y_t|\gs 2+2u)du=2^{-1}\int^{\infty}_{2}\P(|Y_t|\gs v)dv\\
&\ls \int^{\infty}_{2}\frac{v^{2^k-1}}{2^{2^k}}\P(|Y_t|\gs v)dv\ls \frac{1}{2^k N_k}\int^{\infty}_{0}2^k v^{2^k-1}\P(|Y_t|\gs v)dv\\
&=  \frac{\E |Y_t|^{2^k}}{2^k N_k}=\frac{1}{2^kN_k}.
\end{align*}
 Therefore, 
$$
\E\sup_{A\in \ccA_N}|X_{\pi_N(A)}-X_{\pi_0(A)}|\ls 
2\gamma_X(T)(1+\sum^N_{k=1} \frac{N_k}{2^k N_k}) \ls 4\gamma_X(T).
$$
This ends the proof.
\end{dwd}
The question about lower bounds for the suprema of canonical processes is much more involved. Let us summarize the processes in which the full characterization of the supremum (i.e. lower and upper bound) can be provided with the use of $\gamma_X(T)$. The seminal result of Fernique and Talagrand known as the Majorizing Measure Theorem (see \cite{Fer}, \cite{Tal4} or \cite{Tal1} for the modern formulation) is equivalent with the statement that $S_G(T)$ is comparable with $\gamma_G(T)$ up to a numerical constant. In \cite{Tal3} it was proved that $S_X(T)$ is comparable with a quantity which, in a sense, is equivalent to $\gamma_X(T)$ for canonical process generated by $\xi$'s which are symmetric and satisfy $\P(|\xi|>t)=\exp(-c_pt^p)$ for a fixed $p\in [1,2]$. A similar result holds for $p>2$, yet it is only possible to show that there exists a set $T'\subset  \ell^2$ (which may significantly differ from $T$) such that
$S_X(T)$ is comparable with $\gamma_X(T')$ up to a numerical constant.   Note that the limiting case, when $p\ra \infty$ is the question about canonical Bernoulli processes. Later, the idea of \cite{Tal3} was slightly generalized by R. Lata\l{}a in \cite{Lat0} for canonical processes generated by $\xi$ with log-concave tails, yet under specific regularity assumptions. Finally, in \cite{Lat-Tko} it was proved it suffices 
to assume only certain conditions on a moment growth of $\xi$. Unfortunately, this result still does not apply to Bernoulli processes.  The question of characterization of $S_B(T)$ was a long-standing problem posed by M. Talagrand and known as the Bernoulli conjecture. It was finally proved in \cite{Bed1}. In order to explain this result we need to provide a family of distances relevant to canonical Bernoulli processes which follow from some properties of Bernoulli-type random variables. By the results (see \cite{Hit}, \cite{Mon} and \cite{Kwa} for the below formulation) for any $p\in \N$, $p\gs 1$
\be\label{ole0}
\|B_t\|_p\ls \sum^p_{i=1}|t^{\ast}_i|+\sqrt{p}(\sum_{i>p}|t^{\ast}_i|^2)^{\frac{1}{2}}\ls 4\|B_t\|_p,
\ee
where $(t^{\ast}_i)_{i\gs 1}$ is the rearrangement of $(t_i)_{i\gs 1}$ such that
$|t^{\ast}_1|\gs |t^{\ast}_2|\gs \ldots$. Now, if we denote by $I\subset\N$ some index set, we can think of (\ref{ole0}) as a decomposition of the norm $\|B_t\|_p$
into the $\ell^1$ part 
$$
\sum^p_{i=1}|t^{\ast}_i|=\sup_{|I^c|\ls p}\sum_{i\in I^c}|t_i|
$$ 
and the Gaussian part
$$
\sqrt{p}(\sum_{i>p}|t^{\ast}_i|^2)^{\frac{1}{2}}=\sqrt{p}\inf_{|I^c|\ls p}(\sum_{i\in I}|t_i|^2)^{\frac{1}{2}}.
$$
In fact a similar characterization to (\ref{ole0}) can be formulated for a broad class
of processes to mention processes with log-concave distributions. In particular, in \cite{Lat}
there is a characterization of $\|X_t-X_s\|_p$ for canonical processes based on one-unconditional log-concave random variables. As we have mentioned the characterization of $S_B(T)$ was known as the Bernoulli conjecture and was finally proved in \cite{Bed1}. It states that similarly to (\ref{ole0}) the understanding of $S_B(T)$ can be decomposed into the Gaussian and $\ell^1$ part. More precisely, there must exist 
a decomposition of $T$ into $T_1,T_2\subset \ell^2$ such that $T_1+T_2\supset T$ and moreover $S_B(T)$ dominates up to a universal constant both $\sup_{t\in T_1}\|t\|_1$ and $S_G(T_2)$. Usually such decomposition is formulated in the language of existence of 
a mapping $\pi:T\ra \ell^2$ which defines $T_1=\{t-\pi(t):t\in T\}$ and $T_2=\{\pi(t):t\in T\}$. Recall that we can always assume that $0\in T$ and $\pi(0)=0$. We now turn to prove that the Bernoulli Theorem \cite{Bed1} implies that there must exist a subset $T'\subset \ell^2$
such that $\gamma_B(T')$ is comparable to $S_B(T)$. The idea of the proof works also for other classes of canonical 
processes for which we can characterize $S_X(T)$ in terms of increments, see Remark \ref{rema0} below. 
\bt\label{thm:3}
There exists a function $\pi:T\ra \ell^2$ such that
\be\label{atr}
K^{-1}(\gamma_B(T_1+T_2))\ls S_B(T)\ls K(\gamma_B(T_1+T_2)),
\ee
where $K$ is a universal constant, $T_1=\{t-\pi(t):t\in T\}$ and $T_2=\{\pi(t):t\in T\}$.
\et

\begin{dwd}
First, we have to notice that it suffices to prove the result for countable sets $T$. Indeed
for any dense countable set $\bar{T}$ it is true that $S_B(T)=S_B(\bar{T})$. Suppose we have a decomposition of $\bar{T}$
into $\bar{T}_1$ and $\bar{T}_2$ so that (\ref{atr}) holds. It is easy to observe that
$\gamma_B(\bar{T}_1)=\gamma_B(\mathrm{cl}\bar{T}_1)$ and $\gamma_B(\bar{T}_2)=\gamma_B(\mathrm{cl}\bar{T}_2)$ moreover, $T_1$ and $T_2$ must be bounded since otherwise $\gamma_B(T_1)$ or 
$\gamma_B(T_2)$ is infinite and hence also $S_B(T)$. Therefore, $\mathrm{cl}\bar{T}_1+\mathrm{cl}\bar{T}_2$ is compact and contains $\mathrm{cl}T$. Consequently, with no loss in generality we can assume that $T$ is countable.
Then, by the main result of \cite{Bed1} we get the existence of $\pi:T\ra \ell^2$ and consequently the existence of the decomposition into countable sets $T_1,T_2\subset \ell^2$ such that
$T\subset T_1+T_2$ and
\be\label{ole3}
S_B(T)\gs K^{-1}(\sup_{t\in T_1}\|t\|_1+S_G(T_2)),
\ee
where $K$ is a universal constant. By the Pisier's \cite{Pis} and Talagrand's theorems \cite{Tal1} we have that $S_G(T)$ is comparable with $\gamma_G(T)$. Let $g$ be a standard normal variable independent of $B_t$, $t\in T$. Observe that for any $p\gs1$
$$
\frac{\sqrt{2}}{\sqrt{\pi}}\|B_t-B_s\|_{p}= \E |g|\|B_t-B_s\|_p\ls \|G_t-G_s\|_p
$$
and hence $\frac{\sqrt{\pi}}{\sqrt{2}}\gamma_G(T_2)\gs \gamma_{B}(T_2)$. On the other hand,
we can choose an admissible sequence $(T^n_{1})^{\infty}_{n=0}$ such that $\bigcup^{\infty}_{n=0}T^n_1=T_1$.
Fix any given point $t_0$ in $T$. Define 
$$
\pi_m(t)=\left\{\begin{array}{lll}
t_0 & \mbox{if} & t\in T\backslash T^{m} \\
t   & \mbox{if} & t\in T^m\backslash T^{m-1}.
\end{array}\right.
$$ 
If $t\in T^m\backslash T^{m-1}$
$$
\sum^{\infty}_{n=1}\|B_{\pi_n(t)}-B_{\pi_{n-1}(t)}\|_{2^n}=\|B_t-B_{t_0}\|_{2^m}\ls \|t-t_0\|_{1}. 
$$
Therefore, by the triangle inequality
$$
\sup_{t\in T }\sum^{\infty}_{n=1}\|B_{\pi_n(t)}-B_{\pi_{n-1}(t)}\|_{2^n}\ls 2\sup_{t\in T_1}\|t\|_1.
$$
In this way we have proved that
$$
S_B(T)\gs K^{-1}(\gamma_B(T_1)+\gamma_B(T_2))\gs(KC_1)^{-1}\gamma_B(T_1+T_2).
$$
On the other hand, we have a trivial upper bound
$$
S_B(T)\ls S_B(T_1)+S_B(T_2)=S_B(T_1+T_2)\ls 4\gamma_B(T_1+T_2),
$$
by Theorem \ref{theo2}.
\end{dwd}
Let us also observe that for $\P(|\xi_i|>t)=\exp(-c_pt^p)$, $p\gs 2$ we could give a similar proof. 
It is based on the fact that for any $p$ there is a Talagrand's \cite{Tal1} characterization of $S_X(T)$. 
\begin{rema}\label{rema0}
For the class of canonical processes based on independent symmetric $\xi_i$ such that
$\P(|\xi_i|>t)=\exp(-c_pt^p)$, $p\gs 2$, $S_X(T)$ is  comparable with $\gamma_X(T_1+T_2)$ up to a constant for some $T_1+T_2\subset \ell^2$ that contains $T$. The role of $T_2$ may be again addressed to the Gaussian reason, whereas
$T_1\subset \ell^{p^{\ast}}$ for $p^{\ast}=\frac{p}{p-1}$.
\end{rema}
In general, we conjecture that the same is true for canonical processes based on log-concave random variables. 
\begin{conj}

If $(\xi_i)$, $i\gs1$ is a sequence of independent log-concave random variables with mean 0 and variance 1 then there exists $\pi:T\ra \ell^2$
and sets $T_1=\{t-\pi(t)\in \ell^2:\;t\in T\}$, $T_2=\{\pi(t)\in \ell^2:\;t\in T\}$ such that
$$
K^{-1}(\gamma_X(T_1+T_2))\ls S_X(T) \ls K(\gamma_X(T_1+T_2)),
$$
where $K$ is a universal constant.
\end{conj}

\section{Contractions of canonical Bernoulli processes}
Suppose we have a map $\varphi:T\ra\ell^2$. The main question we treat in this paper is under what assumptions on $X_t$, $T$ and $\varphi$ we can show that $S_X(\varphi(T))$ is bounded by $S_X(T)$ up to a numerical constant. In particular we are interested in the case of canonical Bernoulli processes. Let's start with classic results concerning comparison of Gaussian processes. It is well-known that if $G_t$ and $G'_t$, $t\in T$
are centered Gaussian processes and $\E |G_t-G_s|^2\ls \E |G'_t-G'_s|^2$, then for each finite subset $F\subset T$
\be\label{gauss1}
\E \sup_{t\in F} G_t\ls \E \sup_{t\in F} G'_t.  
\ee
This comparison is a consequence of Slepian's Lemma (Corollary 3.14 in \cite{L-T} provides the proof with constant 2, the proof with the best possible constant $1$ is  Corollary 2.1.3 in \cite{Fer}). Note also that by the Majorizing Measure Theorem the result can be generalized to the case when we compare a centered Gaussian process with a centered process for which we only require sub-gaussianity property, see Theorem 12.16 in  \cite{L-T}. We start with a discussion on possible extensions of this result. It is natural to ask for other cases when similar comparison results hold. From Theorem \ref{theo2} it can be easily deduced that if we can compare moments then we can compare $\gamma$-type upper bounds.
\begin{coro}\label{thm:0}
Suppose that $(X_t)_{t\in T}$ is a canonical process and suppose that for each $n\gs1$ , $\varphi:T\ra \ell_2$ and constant $C$ it
satisfies
\be\label{ole1}
\|X_{\varphi(t)}-X_{\varphi(s)}\|_{2^n}\ls C\|X_t-X_s\|_{2^n}, 
\ee
then $S_X(\varphi(T))\ls 4C\gamma_X(T)$.
\end{coro}
\begin{dwd}
Clearly, by Theorem \ref{theo2} we have $S_X(\varphi(T))\ls 4\gamma_X(\varphi(T))\ls 4C\gamma_X(T)$.
\end{dwd}
This means that if we could show that $S_X(t)\gs K^{-1}\gamma_X(T)$, then by Corollary \ref{thm:0}
we would be able to prove that $S_X(\varphi(T))\ls 4CKS_X(T)$. Unfortunately, in general, there is no proof that $\gamma_X(T)$ is comparable with $S_X(T)$. On the other hand, as it was discussed before there are cases where the idea works. In particular, we could use Corollary \ref{thm:0} in order to recover the Gaussian comparison result with some absolute constant. However, in the Gaussian setting, one can simply refer to (\ref{gauss1}) rewriting it in the following way
\be\label{gauss2}
\mbox{if}\;\;\varphi:T\ra \ell^2\;\;\mbox{satisfies}\;\;\|\varphi(t)-\varphi(s)\|_2\ls \|t-s\|_2, 
\mbox{then}\;\;S_G(\varphi(T))\ls S_G(T). 
\ee
\noindent 
We now move to the case of canonical Bernoulli processes. The only known comparison result is Theorem 2.1 in \cite {Tal2} and Theorem 4.12 in \cite{L-T}. It states that if
$\varphi=(\varphi_i)_{i\gs 1}:T\ra \ell_2$, where $\varphi_i:\R\ra \R$ are contractions then $S_B(T)$ dominates $S_B(\varphi(T))$
with the constant $1$, namely
\be\label{bernoulli}
\mbox{if}\;\; |\varphi_i(x)-\varphi_i(y)|\ls |x-y|,\;\; \mbox{for} \;\;i\gs 1\;\;\mbox{then}\;\;S_B(\varphi(T))\ls S_B(T).
\ee
Note that if we are interested in the comparison up to a numerical constant (not necessarily equal $1$) then the requirement of coordinate contractions is too demanding. However, it is known that the result analogous to (\ref{gauss1}), where we assume that $\varphi:\ell^2\rightarrow\ell^2$ is a Lipschitz contraction does not hold for Bernoulli processes. Therefore some additional assumptions on $\varphi$ or $T$ are required. As we show in this paper, the comparison for canonical Bernoulli processes should depend on a suitable family of distances already presented in (\ref{ole0}). The straightforward consequence of Theorem \ref{thm:3} is the following comparison result.
\begin{coro}
Suppose that $\varphi:T\ra\ell^2$ can be extended to $T_1+T_2$ in such a way that for any $p\gs1$
$$
\|B_{\varphi(t)}-B_{\varphi(s)}\|_{p}\ls \|B_t-B_s\|_{p},\;\;\mbox{for all}\;\;s,t\in T_1+T_2
$$
then $S_B(\varphi(T))\ls KS_B(T)$, where $K$ is a universal constant.
\end{coro}
\begin{dwd}
Clearly, by Theorem \ref{theo2} we have $S_B(\varphi(T))\ls 4\gamma_B(\varphi(T))$. Hence, by Theorem \ref{thm:3} 
$$
S_B(\varphi(T))\ls 4\gamma_B(\varphi(T))\ls 4\gamma_B(\varphi(T_1+T_2))\ls 4\gamma_B(T_1+T_2)\ls 4KS_B(T). 
$$
\end{dwd}
Note that the trouble with application of the above result is that  $T_1+T_2$ may be much larger than $T$. We conjecture the following generalization of the above result.

\begin{conj}
Let $\varphi=(\varphi_i)_{i\gs 1}:T\ra \ell^2$. If  
\be\label{ole101}
\|B_{\varphi(t)}-B_{\varphi(s)}\|_p\ls \|B_t-B_s\|_{p},\;\;\mbox{for all}\;\;p\gs 2,\;s,t\in T
\ee
then $S_B(\varphi(T))\ls KS_B(T)$, for an absolute constant $K$.
\end{conj}
Towards this aim we prove a weaker form of the conjecture. 
As we have explained the norm $\|B_t-B_s\|_p$ can be decomposed into the Gaussian and $\ell^1$ part.
Our condition states that if Gaussian part of $\|B_t-B_s\|_p$ dominates Gaussian part of $\|B_{\varphi(t)}-B
_{\varphi(s)}\|_p$, 
for all $s,t\in T$ and  $p\gs 1$ then $S_B(T)$ dominates $S_B(\varphi(T))$ up to an absolute constant.
\bt\label{thm:1}
Suppose that for all $s,t\in T$ and all natural $p$ such that $p\gs 0$  we have
\be\label{ole2}
\inf_{|I^c|\ls Cp}\sum_{i\in I}|\varphi_i(t)-\varphi_i(s)|^2\ls C^2\inf_{|I^c|\ls p}\sum_{i\in I}|t_i-s_i|^2,
\ee
for an absolute constant $C\gs 1$. Then $S_B(\varphi(T))\ls KS_B(T)$, where $K$ is a universal constant.
\et
The result is stronger than the comparison for Bernoulli processes (\ref{bernoulli}).
In this way Theorem \ref{thm:1} supports the conjecture that (\ref{ole101}) suffices to prove that $S_B(\varphi(T))\ls KS_B(T)$.
Note that there is an important  case for which the conjecture is true. Namely, when we assume that all supports 
$J(t)=\{i\gs 1:\;|t_i|>0\}$ of $t\in T$ are disjoint. It is crucial is to understand that in this case the decomposition 
postulated in the Bernoulli Theorem can have a special form: 
$\pi(t)=t_{J^1(t)}$ and $t-\pi(t)=t_{J^2(t)}$, where $J^1(t)$ and $J^2(t)$
are disjoint and $J^1(t)\cup J^2(t)=J(t)$. We show this fact
when proving  the following result.
\bt\label{thm:5}
Suppose that (\ref{ole101}) is satisfied and supports $J(t)=\{i\gs 1:\;|t_i| >0\}$ are disjoint for all $t\in T$ then
$S_B(\varphi(T))\ls KS_B(T)$, where $K$ is a universal constant.
\et
As we show in the last section, results of this type are of interest when one wants to compare weak and strong moments for random series in a Banach space.
The question was proposed by K. Oleszkiewicz in private communication.

\section{Proof of the main result}

In this section we prove Theorem \ref{thm:1} and Theorem \ref{thm:5}.  
\begin{dwd}[Proof of Theroem \ref{thm:1}]
The main step in the proof of the Bernoulli theorem - Proposition 6.2  in \cite{Bed1} is to show the existence of a suitable admissible
sequence of partitions. Consequently, if $S_B(T)<\infty$ and $0\in T$ then it is possible to define
nested partitions $\ccA_n$ of $T$ such that $|\ccA_n|\ls N_n$. Moreover, for each $A\in \ccA_n$ it is possible to find $j_n(A)\in \Z$ and $\pi_n(A)\in T$ (we use the 
notation $j_n(t)=j_n(A_n(t))$ and $\pi_n(t)=\pi_n(A_n(t))$, where $t\in A_n(t)\in \ccA_n$) which satisfy
the following conditions
\begin{enumerate}[(i)]
\item $\|t-s\|_2\ls \sqrt{M}r^{-j_0(T)}$, for $s,t\in T$;
\item if $n\gs 1$, $\ccA_n\ni A\subset A'\in \ccA_{n-1}$ then 
\begin{enumerate}
\item either $j_n(A)=j_{n-1}(A')$ and $\pi_n(A)=\pi_{n-1}(A')$ 
\item or $j_n(A)>j_{n-1}(A')$, $\pi_n(A)\in A'$ and
\be\label{omyk}
\sum_{i\in I_n(A)}\min\{|t_i-\pi_n(A)_i|^2,r^{-2j_n(A)}\}\ls M 2^nr^{-2j_n(A)},
\ee
where for any $t\in A$
$$
I_n(A)=I_n(t)=\{i\gs 1:\;|\pi_{k+1}(t)_i-\pi_k(t)_i|\ls r^{-j_k(t)}\;\;\mbox{for}\;\;0\ls k\ls n-1\}.
$$
\end{enumerate}
\item Moreover, numbers $j_n(A)$, $A\in \ccA_n$, $n\gs 0$ satisfy
\be\label{ole4}
\sup_{t\in T}\sum^{\infty}_{n=0}2^nr^{-j_n(t)}\ls LS_B(T),
\ee
where $L$ is an absolute constant.
\end{enumerate}
As proved in Theorem 3.1 in \cite{Bed1} the existence of the quantities $\ccA_n,j_n(A),\pi_n(A),I_n(A)$ that satisfy conditions (\emph{i}) and (\emph{ii}) formulated above implies the existence of a decomposition $T_1,T_2\subset \ell^2$, $T_1+T_2\supset T$ such that
$$
\sup_{t^1\in T_1}\|t^1\|_1\ls LM\sup_{t\in T}\sum^{\infty}_{n=0}2^nr^{-j_n(t)}\;\;\mbox{and}\;\;
\gamma_G(T_2)\ls L\sqrt{M}\sup_{t\in T}\sum^{\infty}_{n=0}2^nr^{-j_n(t)}.
$$    
Together with the condition (\emph{iii}) we get (\ref{ole3}). Our aim is to use the mapping $\varphi$ to transport all the required quantities
to $\varphi(T)$. Before we do it we formulate an auxiliary fact about sets $I_n(A)$, namely we show that we can get rid of truncation
in (\ref{omyk}) if we skip a well controlled number of coordinates. 
We observe that for each $t\in A\in \ccA_n$ there must exist set $J_n(t)$ such that $|J^c_n(t)|\ls M 2^{n+1}$
and 
\be\label{omsk}
\sum_{i\in J_n(t)}|t_i-\pi_n(t)_i|^2\ls M 2^nr^{-2j_n(t)}.
\ee
The fact will be proved in two steps. First, we show that $|I_n(t)^c|\ls M2^{n}$.   
We may only prove that $|I_n(t)|=|I_n(A_n(t))|\ls 2^{n}$, if $\pi_{n-1}(t)\neq \pi_n(t)$, which
implies $j_{n-1}(t)\neq j_n(t)$ and $\pi_n(t)\in A_{n-1}(t)$. Therefore, there exists $k\in \{1,\ldots,n\}$
such that
$$
j_{n-1}(t)=j_{n-k}(t)>j_{n-k-1}(t),\;\;\mbox{where use the notation}\;\;j_{-1}(t)=-\infty 
$$  
and hence $\pi_{n}(t)\in A_{n-1}(t)\subset A_{n-k}(t)$ and $\pi_{n-1}(t)=\pi_{n-k}(t)$, $j_{n-1}(t)=j_{n-k}(t)$ so by the
construction of $(\ccA_n)_{n\gs 0}$
\begin{align*}
& \sum_{i\in I_{n-k}(t)}\min\{(\pi_{n}(t)_i-\pi_{n-1}(t)_i)^2,r^{-2j_{n-1}(t)} \}\\
&=\sum_{i\in I_{n-k}(t)}\min\{(\pi_{n}(t)_i-\pi_{n-k}(t)_i)^2,r^{-2j_{n-k}(t)}\}\ls M2^{n-k}r^{-2j_{n-k}(t)}.
\end{align*}
Consequently,
$$
|\{i\in I_{n-k}(t):\;|\pi_{n}(t)_i-\pi_{n-1}(t)_i|>r^{-j_{n-1}(t)}\}|\ls M 2^{n-k}. 
$$
Obviously, 
$$
|I_n(t)^c|\ls |I_{n-k}(t)^c|+|\{i\in I_{n-k}(t):\;|\pi_{n}(t)_i-\pi_{n-1}(t)_i|>r^{-j_{n-1}(t)}\}|\ls |I_{n-k}(t)^c|+M 2^{n-k}.
$$
Therefore by the induction, $|I_n(t)^c|\ls M\sum^n_{k=1}2^{n-k}\ls M 2^{n}$. Let
$$
J_n(t)=\{i\in I_n(A):\; |t_i-\pi_n(t)_i|\ls r^{-j_n(A)} \}.
$$
The second step is to establish that $|I_n(t)\backslash J_n(t)|\ls M2^{n}$. 
Again it suffices to prove the result only for $n$ such that $j_n(t)>j_{n-1}(t)$. 
Note that by (\ref{omyk})
$$
|I_n(t)\backslash J_n(t)|r^{-2j_n(t)}=\sum_{i\in I_n(A)\backslash J_n(t)}r^{-2j_n(t)}\ls M2^nr^{-2j_n(t)}
$$
and hence the result holds. It remains to observe that 
$$
|J_n(t)^c|\ls |I_n(t)^c|+|I_n(t)\backslash J_n(t)|\ls M(2^{n}+2^n)\ls M2^{n+1}.
$$
We turn to construct an admissible partition sequence together with all the supporting quantities for the set $\varphi(T)$. 
Let $\ccB_{n}$ consists of $\varphi(A)$, $A\in \ccA_n$. Obviously partitions $\ccB_n$ are admissible, nested
and $\ccB_0=\{\varphi(T)\}$. Moreover, for each $n\gs 0$ and $A\in \ccA_n$ we define
$$
\pi_n(\varphi(A))=\varphi(\pi_n(A))\;\;\mbox{and}\;\;j_n(\varphi(A))=j_n(A)
$$ 
and obviously
\begin{align*}
& I_n(\varphi(A))=I_n(\varphi(t))\\
&=\{i\gs 1: |\varphi(\pi_{k+1}(t))_i-\varphi(\pi_k(t))_i|\ls r^{-j_k(\varphi(t))}\;\mbox{for}\;0\ls k\ls n-1\}.
\end{align*}
As we have mentioned at the beginning of this proof, in order to use Theorem 3.1 in \cite{Bed1} we have to verify conditions (\emph{i}) and (\emph{ii}) for the new sequence $\ccB=(\ccB_n)_{n\gs 0}$
as well as $j_n(B),\pi_n(B),I_n(B)$ for $B\in \ccB_n$, $n\gs 0$.
For this aim we need our main condition (\ref{ole2}). 
First it is obvious that  that (\ref{ole2}) implies for $p=0$ that
$$
\|\varphi(t)-\varphi(s)\|_2\ls \|t-s\|_2\ls\sqrt{M}r^{-j_0(T)}.
$$
If $A\in\ccB_n$ and $\varphi(A)\subset \varphi(A')\in \ccB_{n-1}$ then either
$$
j_n(\varphi(A))=j_n(A)=j_{n-1}(A')=j_{n-1}(\varphi(A'))
$$
and
$$
\pi_n(\varphi(A))=\varphi(\pi_n(A))=
\varphi(\pi_{n-1}(A'))=\pi_{n-1}(\varphi(A'))
$$
or $j_n(\varphi(A))=j_n(A)>j_{n-1}(A')=j_{n-1}(\varphi(A'))$. In this case we have
$\pi_n(\varphi(A))=\varphi(\pi_n(A))\in \varphi(A')$ and it suffices to show that
\be\label{oleg5}
\sum_{i\in I_n(\varphi(A))}\min\{|\varphi(t)_i-\varphi(\pi_n(A))_i|^2,r^{-2j_n(\varphi(A))}\}\ls C2^nr^{-2j_n(\varphi(A))}.
\ee
Obviously, the problem now is that we know a little about the structure of the set $I_n(\varphi(A))$.
Therefore, we simply prove that 
$$
\sum_{i\gs 1}\min\{|\varphi(t)_i-\varphi(\pi_n(A))_i|^2,r^{-2j_n(\varphi(A))}\}\ls C2^nr^{-2j_n(\varphi(A))}.
$$
It is obvious that 
\begin{equation}\label{omsk3}
\sum_{i\gs 1}\min\{|\varphi(t)_i-\varphi(\pi_n(A))_i|^2,r^{-2j_n(\varphi(A))}\}\\
\ls C_2 2^n r^{-2j_n(A)}+\inf_{|I^c|\ls C_2 2^n}\sum_{i\in I}|\varphi(t)_i-\varphi(\pi_n(A))_i|^2.
\end{equation}
We can choose $C_2\gs 2CM$ in a way that by (\ref{ole2}) we get
\begin{align*}
& \inf_{|I^c|\ls C_2 2^n}\sum_{i\in I}|\varphi(t)_i-\varphi(\pi_n(A))_i|^2\\
&\ls C^2\inf_{|I^c|\ls M2^{n+1}}\sum_{i\in I}|t_i-\pi_n(A)_i|^2\ls C^2\sum_{i\in J_n(t)}|t_i-\pi_n(A)_i|^2.
\end{align*}
Hence, by (\ref{omsk}) and (\ref{omsk3})
$$
\sum_{i\gs 1}\min\{|\varphi(t)_i-\varphi(\pi_n(A))_i|^2,r^{-2j_n(\varphi(A))}\}\ls (C_2+C^2 M) 2^n r^{-2j_n(A)}
$$ 
which proves (\ref{oleg5}) with $C_3=C_2+C^2 M$. We have proved that assumptions required in Theorem 3.1 in \cite{Bed1} are satisfied for $(\ccB_n)_{n\gs 0}$
and the supporting quantities. Consequently, there exists a decomposition $S_1,S_2\subset \ell^2$ such that $S_1+S_2\supset \varphi(T)$
and
$$
\sup_{s\in S_1}\|s\|_1\ls LC\sup_{t\in \varphi(T)}\sum_{n\gs 0}2^nr^{-j_n(t)},\;\;
\gamma_G(S_2)\ls L\sqrt{C}\sup_{t\in \varphi(T)}\sum_{n\gs 0}2^nr^{-j_n(t)}. 
$$
Since $j_n(\varphi(t))=j_n(t)$ and we have (\ref{ole4}) for $(\ccA_n)_{n\gs 0}$ we obtain that
$$
\sup_{t\in \varphi(T)}\sum_{n\gs 0}2^nr^{-j_n(t)}\ls LS_B(T).
$$
It implies that
$$
S_B(\varphi(T))\ls S_B(S_1)+S_B(S_2)\ls KS_B(T),
$$
for a universal constant $K$ and ends the proof.
\end{dwd}
The second case we consider is when for all $t\in T$ supports $J(t)=\{i\gs 1:\; |t_i|>0\}$ are disjoint. The proof requires the following notation. For any $t\in \ell^2$
and $J\subset \{1,2,\ldots\}$ we define $t1_J\in \ell^2$ such that $(t1_J)_i=t_i$ for $i\in J$ and $(t1_J)_i=0$ otherwise.

\begin{dwd}[Proof of Theorem \ref{thm:5}]
Obviously, we may require that $b(T)<\infty$. We additionally assume that $0\in T$. It simplifies the proof, but it works also for the general case as we will point out at the end. Recall that by Bernoulli Theorem \cite{Bed1}
there exists a decomposition $T_1+T_2\supset T$ such that
\begin{align}\label{bound}
S_B(T)\gs K^{-1}(\sup_{t\in T_1}\|t\|_1+\gamma_G(T_2)),
\end{align}
where $K$ is an absolute constant. Obviously, we may think of $K$
as suitably large. We can represent the decomposition by $\pi:T\ra \ell^2$ in a way that  
$T_2=\{\pi(t):\;t\in T\}$ and $T_1=\{t-\pi(t):\;t\in T\}$. 
We show that under the disjoint supports assumption we
may additionally require
that $\pi(t)=t1_{J^2(t)}$ and $t-\pi(t)=t1_{J^1(t)}$ where $J^1(t)$ and $J^2(t)$
are disjoint subsets of $J(t)$ such that $J^1(t)\cup J^2(t)=J(t)$.
Moreover, $J^2(t)=\{i\in J(t):\;|t_i|\ls p(t)\}$, for some suitably chosen $p(t)\gs 0$.
\smallskip

\noindent In order to prove the result we have to look closer into 
the definition of $\pi(t)$ in the proof of Theorem 3.1 in \cite{Bed1}.
The definition is based on the construction of admissible partitions
we have described in the proof of Theorem \ref{thm:1} above. Using the notation introduced there let
\be\label{humbak}
m(t,i)=\inf\{n\gs 0:\;\;|\pi_{n+1}(t)_i-\pi_n(t)_i|>r^{-j_n(t)}\},\;\;t\in T,i\gs 1.
\ee
Note that $S_B(T)$ is comparable with $\sup_{t\in T}\sum_{n\gs 0}2^nr^{-j_n(t)}$. Therefore, if $S_B(T)$ is finite then necessarily $\lim_{n\ra\infty}j_n(t)=\infty$ for all $t\in T$. From the partition construction used in Section 6 in \cite{Bed1} we know that we can additionally assume a regularity condition on $j_n(t)$, $n\gs0$, namely
$$j_n(t)\ls j_{n-1}(t)+2 \;\; \mbox{for all}\;\; n\gs0$$ 
and for technical purpose we take $j_{-1}(t)=-\infty$.
As in the proof of Theorem 3.1 in \cite{Bed1} the Bernoulli decomposition $\pi(t)$ is given by 
$\pi(t)_i=\pi_{m(t,i)}(t)_i$, where if $m(t,i)=\infty$ the definition means that $\pi(t)_i=\lim_{n\ra\infty}\pi_n(t)_i$
and the limit exists.  Consequently, denoting $J_n(t)=\{i\gs 1:\;m(t,i)=n\}$ and $J_{\infty}(t)=\{i\gs 1:\;m(t,i)=\infty\}$ we get
$$
\pi(t)=\sum_{n\gs 0}\pi_n(t)1_{J_n(t)}+\pi(t)1_{J_{\infty}(t)}.
$$
Clearly, $J_n(t)$, $n\gs 0$ and $J_{\infty}(t)$ are disjoint.  Note also that if $m(t,i)=\infty$ and $i\in J(\pi(t))$, then
there must exist $n\gs 0$ such that $|\pi_k(t)_i|>0$ for all $k\gs n$. Due to the disjoint supports assumption 
it is only possible if there exists $n\gs 0$ such that $\pi_n(t)_i=\pi_{n+1}(t)_i=\ldots$. Now, 
if there exists $m\gs 0$ such that $A_m(t)=\{t\}$ we define
$$
\tau(t)=\inf\{n\gs 0:\; A_n(t)=\{t\}=\{\pi_n(t)\},\;j_{n-1}(t)<j_n(t)\},\mbox{otherwise}\;\;\tau(t)=\infty.
$$
The moment $\tau(t)$ is of special nature in the sense that without loss of generality we may assume that for $n\gs\tau(t)$ it is true that $j_n(t)=j_{n-1}(t)+2$. It is due to the fact the partition is ceased after this moment.
Now, we define
$$
J^2(t)=\{i\in J(t):\; |t_i|\ls r^{-j_{\tau(t)}(t)-1} \},\;\;J^1(t)=J(t)\backslash J^2(t).
$$
We can now introduce the improved version of $\pi$ denoted by $\bar{\pi}$ and given by
$$
\bar{\pi}(t)=t 1_{J^2(t)}.
$$  
It is clear that 
$$
\|t-\bar{\pi}(t)\|_1=\|t 1_{J^1(t)}\|_1.
$$
For $n\gs 0$ let 
$$
L_n(t)=\{i\in J(t):\; r^{-j_n(t)}<|t_i|\ls r^{-j_{n-1}(t)}\}.
$$
Observe, that  $J^1(t)= \bigcup_{n< \tau(t)} L_n(t)$.
If $i\in L_n(t)$, $n\gs 0$, then 
we may find $0\ls m\ls n$ such that $j_{m-1}(t)<j_m(t)=j_{m+1}(t)=\ldots =j_n(t)$. 
Consequently, using the definition (\ref{omyk}) of $I_n(t) $  for all $s\in A_m(t)$ 
\begin{align*}
& \sum_{i\in I_{n}(t)}\min\{|s_i-\pi_n(t)_i|^2,r^{-2j_n(t)}\}=\sum_{i\in I_{m}(t)}\min\{|s_i-\pi_m(t)_i|^2,r^{-2j_m(t)}\}\\
&\ls M2^mr^{-2j_m(t)}=M2^m r^{-2j_n(t)}\ls M2^n r^{-2j_n(t)}.
\end{align*}
We need to show that the decomposition $\bar{\pi}$ is of the right form i.e. satisfies (\ref{bound}). For this aim we need to investigate a few cases following from different possible paths of approximations $\pi$.
First suppose that $t\neq \pi_n(t)$. Then we may use the above inequality for $s=t$ and due to the disjoint supports we have 
$$
|I_n(t)\cap L_n(t)|r^{-2j_n(t)}\ls
\sum_{i\in I_{n}(t)}\min\{|s_i-\pi_n(t)_i|^2,r^{-2j_n(t)}\}\ls
 M2^nr^{-2j_n(t)},\;\;\mbox{so}\;\;|I_n(t)\cap L_n(t)|\ls M2^n.
$$
The same inequality holds if $t=\pi_n(t)$ but $A_m(t)\neq \{t\}$. 
We show that $L_n(t)\subset I_n(t)$. Indeed, suppose that $i\not\in I_n(t)$. It means that for some $k\in \{0,1,\ldots, n-1\}$
we have $|\pi_{k+1}(t)_i-\pi_{k}(t)_i|>r^{-j_k(t)}$. This may concern $i\in J(t)$ only if $\pi_{k+1}(t)=t,\pi_k(t)\neq t$
or $\pi_k(t)=t$ and $\pi_{k+1}(t)\neq t$, but then it means that $|t_i|>r^{-j_k(t)}\gs r^{-j_{n-1}(t)}$ i.e. $i\not\in L_n(t)$. It concludes the argument that $L_n(t)\subset I_n(t)$. 
For $1\ls n<\tau(t)$ it implies that 
\be\label{ppp1}
\sum_{i\in L_n(t)}|t_i|\ls M2^nr^{-j_{n-1}(t)}.
\ee
For $n=0$ we use simply that $|t_i|\ls 2S_B(T)$ and hence 
\be\label{ppp2}
\sum_{i\in L_0(t)}|t_i|\ls 2MS_B(T).
\ee
Now suppose that $t=\pi_n(t)=\pi_m(t)$ and $A_m(t)=\{t\}$. If either $t\neq \pi_{m-1}(t)$
or $\{t\}\neq A_{m-1}(t)$, then $\tau(t)=m$. Otherwise $\tau(t)<m$.
If $\tau(t)=m$, then by the above argument
 $$
\sum_{i\in L_n(t)}|t_i|^2=\sum_{i\in L_n(t)}\min\{|t_i|^2,r^{-2j_{m-1}(t)}\} \ls
 M2^{m-1}r^{-2j_{m-1}(t)},
$$
and thus using that $|t_i|\gs r^{-j_m(t)-1}$ and $j_{m}(t)= j_{m-1}(t)+2$, we have
$$
\sum_{i\in L_n(t)}|t_i|\ls M2^{m-1}r^{-2j_{m-1}(t)+j_m(t)}\ls M2^{m-1}r^{-j_{m-1}(t)+2}.
$$
We have the remaining bound
\be\label{ppp3}
\sum_{i\in L_n(t)}|t_i|\ls M2^{\tau(t)-1}r^{-j_{\tau(t)-1}(t)+2}.
\ee
Combining (\ref{ppp1}), (\ref{ppp2}) and (\ref{ppp3}) we conclude by (\ref{ole4})
\be\label{omar1}
\|t1_{J^1(t)}\|_1\ls 2MS_B(T)+2M\sum^{\tau(t)-2}_{n=0}r^{-j_n(t)}2^n
+M2^{\tau(t)-1}r^{-j_{\tau(t)-1}(t)+2} \ls 2MLS_B(T),
\ee
where $L$ is an absolute constant.
\smallskip

\noindent
Now consider $s,t\in T$, $s\neq t$.
In order to prove that
\be\label{omar2}
\|\bar{\pi}(s)-\bar{\pi}(t)\|_{2}=\|t1_{J^2(t)}-s1_{J^2(s)}\|_2\ls \|\pi(t)-\pi(s)\|_2
\ee
we have to argue that $J^2(t)\cap J(\pi(s))=\emptyset$, $J^2(s)\cap J(\pi(t))=\emptyset$ 
for all $n\gs 0$. Note that $J^2(t)\subset J_{\infty}(t)$ and $J^2(s)\subset J_{\infty}(s)$. Moreover, $J_{\infty}(s)$ and $J_{\infty}(s)$ are disjoint. Obviously, it suffices to show the argument that $J^2(t)\cap J(\pi(s))=\emptyset$.
\smallskip

\noindent
First, note that $J^2(t)\cap J_{\infty}(s)=\emptyset$. Indeed
if the set was non-empty then for a given $n\gs 0$ we would have $t=\pi_n(s)=\pi_{n+1}(s)=\ldots$, but then $s\in A_n(t)$ for all $n\gs 0$ and therefore $\tau(t)=\infty$. This would imply $J^2(t)=\emptyset$ which is a contradiction. Suppose that $i\in J^2(t)$ and $i\in J_n(s)$. This is only possible if $\pi_n(s)=t$ and $\pi_{n+1}(s)\neq \pi_n(s)=t$ and $r^{-j_{n}(s)}<|\pi_n(s)_i|$. Let $m\gs 0$ be such that $j_{m-1}(s)<j_m(s)=j_{m+1}(s)=\ldots=j_n(s)$, then either $m=0$ or $m\gs 1$ and
$t=\pi_n(s)=\pi_m(s)\in A_{m-1}(s)$, which means that $A_{m-1}(s)=A_{m-1}(t)$ and $j_{m-1}(s)=j_{m-1}(t)$. Therefore, $\tau(t)\gs m$ and $j_{\tau(t)}(t)>j_{m-1}(t)$. If $i\in J^2(t)\cap J_n(s)$, then 
$$
r^{-j_{m-1}(t)-2}=r^{-j_{m-1}(s)-2}\ls r^{-j_{m}(s)}<|t_i|\ls r^{-j_{\tau(t)}(t)}=r^{-j_{\tau(t)-1}(t)-2}\ls r^{-j_{m-1}(t)-2},
$$
which is a contradiction. If $m=0$, then the argument is trivial. 
\smallskip

\noindent
Summing up, by (\ref{omar1}) we have
$$
\sup_{t\in T}\|t-\bar{\pi}(t)\|_1\ls LS_B(T)
$$ 
and  by (\ref{omar2}) and the Gaussian comparison we have
$\gamma_G(\pi(T))\ls \gamma_G(\bar{\pi}(T))$, which means that our improved 
version of $\pi$ satisfies
$$ 
S_B(T)\gs K^{-1}(\sup_{t\in T}\|t-\bar{\pi}(t)\|_1+\gamma_G(\bar{\pi}(T))),
$$
where $K$ is a universal constant. In this way we have proved that we may additionally 
require that $\pi(t)=t1_{J^2(t)}$ and $t-\pi(t)=t1_{J^1(t)}$ for some disjoint $J^1(t),J^2(t)$
such that $J^1(t)\cup J^2(t)=J(t)$.  Recall that
$J^2(t)$ in each case is of the form $\{i\in J(t):\; |t_i|\ls r(t)\}$, for a given $r(t)\gs 0$.
\smallskip

\noindent
We turn to the main part of the proof. 
Let $p(t)$ be the smallest positive integer such that
\be\label{omar3}
\sqrt{p(t)}\|t1_{J^2(t)}\|_2\gs KS_B(T)\gs \|t1_{J^1(t)}\|_1.
\ee
Note that it is possible that $J^2(t)=\emptyset$ in which case we may think of $p(t)$
as equal $\infty$. Since $K$ is large enough and $S_B(T)\gs \frac{1}{2}\sup_{t\in T}\|t\|_2$
it is clear that $p(t)$ must be at least greater than, say, $2$.
Consequently, by the choice of $p(t)$ 
\be\label{omar4}
\sqrt{p(t)}\|t1_{J^2(t)}\|_2\ls  2KS_B(T).
\ee
The last step is to define a suitable decomposition for $\varphi(T)$.
For each $t\in T$ we define $\pi(\varphi(t))=t_{J^2(\varphi(t))}$
and $\varphi(t)-\pi(\varphi(t))=t_{J^1(\varphi(t))}$, where
$J^2(\varphi(t))$ and $J^1(\varphi(t))$ are defined by the decomposition
of the norm $\|B_{\varphi(t)}\|_{p(t)}$ i.e.
$$
\sum_{i\in J^1(\varphi(t))}|\varphi(t)_i|=\sup_{|I^c|\ls p(t)}\sum_{i\in I^c}|\varphi(t)_i|
$$
and
$$
\sum_{i\in J^2(\varphi(t))}|\varphi(t)_i|^2=\inf_{|I^c|\ls p(t)}\sum_{i\in I}|\varphi(t)_i|^2.
$$
Consequently by the decomposition (\ref{ole0})  and the main assumption (\ref{ole101}),
\begin{align*}
& \sum_{i\in J^1(\varphi(t))}|\varphi(t)_i|+\sqrt{p(t)}(\sum_{i\in J^2(\varphi(t))}|\varphi(t)_i|^2)^{\frac{1}{2}}\\
&\ls 4\|B_{\varphi(t)}\|_{p(t)}\ls 4\|B_{t}\|_{p(t)}\ls 4(\|t1_{J^1(t)}\|_1+\sqrt{p(t)}\|t1_{J^2(t)}\|_2). 
\end{align*}
Therefore, using (\ref{omar3}), (\ref{omar4})
$$
\sum_{i\in J^1(\varphi(t))}|\varphi(t)_i|\ls K_1S_B(T). 
$$
Moreover, by (\ref{omar3})
$$
(\sum_{i\in J^2(\varphi(t))}|\varphi(t)_i|^2)^{\frac{1}{2}}\ls K_2\|t1_{J^2(t)}\|_2.
$$
It implies that
\begin{align*}
& \|\pi(\varphi(t))-\pi(\varphi(s))\|_2\ls \|\pi(\varphi(t))\|_2+\|\pi(\varphi(s))\|_2\\
&\ls K_2(\|t1_{J^2(t)}\|_2+\|s1_{J^2(s)}\|_2 )\ls K_3\|\pi(t)-\pi(s)\|_2.
\end{align*}
Therefore, by the Gaussian comparison, we get 
$\gamma_G(\pi(\varphi(T)))\ls K\gamma_G(\pi(T))$ and hence finally
$$
S_B(\varphi(T))\ls K(\sup_{t\in T}\|\pi(\varphi(t))\|_1+\gamma_G(\pi(\varphi(T))))\ls KLS_B(T).
$$
It ends the proof in the case when $0\in T$. For the general case the proof follows the same lines, where instead of $t$ we consider $t-\pi_{0}(t)$. Notice that formally this may not obey the disjoint supports assumption, but it does not affect qualitatively the argument presented above.
\end{dwd}
Note that the above proof works since in the case of disjoint supports we have almost perfect knowledge about the decomposition in Bernoulli Theorem. On the other hand, it is not difficult to give an alternative proof based on the independence of variables $B_t$, $t\in T$, but it is worth seeing what the decomposition in Theorem 3.1 in \cite{Bed1} should be in order to make Bernoulli comparison possible.

\section{The Oleszkiewicz problem}

In this section we give an example how to apply our result to compare
expectations of norms of random series in a Banach space. First, we prove a general result which concerns
$\varphi:T\ra \ell^2$ where $\varphi$ is linear, $T$ is convex and $T=-T$. Then,
the assumption (\ref{ole1}) becomes 
\be\label{ole7}
\|B_{\varphi(u)}\|_p\ls C\|B_u\|_p\;\;\mbox{for all}\;p\gs 1\;\;\mbox{and}\;\;u\in \mathrm{cl}(\mathrm{Lin}(T)),
\ee
where $\mathrm{Lin}(T)$ is the linear space spanned by the set $T$. It is because by the assumptions on $T$ any point $u\in \mathrm{Lin}(T)$ can be represented 
as $c\cdot t$, where $c\in \R$ and $t\in T$. By the linearity of $\varphi$
$$
\|B_{\varphi(u)}\|_p=|c|\|B_{\varphi(t)}\|_p\ls C|c|\|B_t\|_p=C\|B_{u}\|_p.
$$
On the other hand, we can easily extend the condition (\ref{ole7}) on
the closure of $\mathrm{Lin}(T)$.
We turn to prove that if $\mathrm{cl}(\mathrm{Lin}(T))=\ell^2$ then (\ref{ole7}) implies that $S_B(T)$ dominates $S_B(\varphi(T))$.
\begin{theo}\label{thm:2}
Suppose that $T=-T$, $T$ is convex and $\mathrm{cl}(\mathrm{Lin}(T))=\ell^2$, if $\varphi$ is linear and satisfies (\ref{ole1}) then
$S_B(\varphi(T))\ls KS_B(T)$, where $K$ is a universal constant.
\end{theo}
\begin{dwd}
By the Bernoulli theorem \cite{Bed1} we have that there exist $T_1,T_2$ such that $T\subset T_1+T_2$ and
$$
S_B(T)\gs L^{-1}(\sup_{t\in T_1}\|t\|_1+\gamma_G(T_2)).
$$
Since $\varphi$ is linear it can be easily extended to $\mathrm{cl}(\mathrm{Lin}(T))=\ell^2$ and thus we can define 
$S_i=\varphi(T_i)$, $i\in\{1,2\}$. Obviously $S_1+S_2\supset \varphi(T)$ moreover (\ref{ole7})
implies in particular that
$$
\|\varphi(u)\|_1=\|B_{\varphi(u)}\|_{\infty}\ls C\|B_u\|_{\infty}=C\|u\|_1.
$$
and
$$
\|\varphi(u)-\varphi(v)\|_2=\|B_{\varphi(u-v)}\|_{2}\ls C\|B_{u-v}\|_2=C\|u-v\|_2.
$$
Consequently
$$
\sup_{s\in S_1}\|s\|_1=\sup_{t\in T_1}\|\varphi(t)\|_1\ls C\sup_{t\in T_1}\|t\|_1
$$
and
$$
\gamma_G(S_2)=\gamma_G(\varphi(T_2))\ls C\gamma_G(T_2).
$$
Therefore
\begin{align*}
& S_B(\varphi(T))\ls S_B(S_1)+S_B(S_2)\ls K(\sup_{s\in S_1}\|s\|_1+\gamma_G(S_2))\\
&\ls CK(\sup_{t\in T_1}\|t\|_1+\gamma_G(T_2))\ls CK^2 S_B(T).
\end{align*}
It ends the proof.
\end{dwd}
We aim to study the question posed by Oleszkiewicz that concerns comparability of weak and strong moments for Bernoulli series in a Banach space.  
Let $x_i$, $y_i$, $i\gs 1$ be vectors in  a Banach space $(B,\|\cdot \|)$. 
Suppose that for all $x^{\ast}\in B^{\ast}$ and $u\gs 0$
\be\label{Ole1}
\P(|\sum_{i\gs 1}x^{\ast}(x_i)\va_i|>u)\ls \bar{C}\P(|\sum_{i\gs 1}x^{\ast}(y_i)\va_i|>\bar{C}^{-1}u).
\ee
This property is called weak tail domination. As we have explained
in the introduction the weak tail domination can be understood 
in terms of comparability of weak moments, i.e. for any integer $p\gs 1$ and $x^{\ast}\in B^{\ast}$
\be\label{Ole2}
\| \sum_{i\gs 1}x^{\ast}(x_i)\va_i\|_p \ls C\| \sum_{i\gs 1}x^{\ast}(y_i)\va_i \|_p
\ee  
Oleszkiewicz asked whether or not it implies the comparability of strong moments.
Namely whether (\ref{Ole1}) or rather (\ref{Ole2})  implies that
\begin{align} 
& \E \| \sum_{i\gs 1}x_i\va_i\|=\E \sup_{x^{\ast} \in B^{\ast}_1} \sum_{i\gs 1} x^{\ast}(x_i)\va_i \nonumber\\
\label{Ole0}&\ls K\E \sup_{x^{\ast} \in B^{\ast}_1} \sum_{i\gs 1} x^{\ast}(y_i)\va_i= K \E \|\sum_{i\gs 1} y_i\va_i\|,
\end{align}
where $K$ is an absolute constant.
Note that in the Oleszkiewicz problem one may assume that $B$ is a separable space since we can easily restrict $B$ to the closure of $\mathbf{Lin}(y_1,x_1,y_2,x_2,\ldots)$. Therefore 
we have that
$$
\E \|\sum_{i\gs 1} y_i\va_i\|=\sup_{F\subset B^{\ast}_1}\E\sup_{x^{\ast}\in F}|\sum_{i\gs 1}x^{\ast}(y_i)\va_i|,
$$
where the supremum is taken over all finite sets $F$ contained in $B^{\ast}_1=\{x^{\ast}\in B^{\ast}:\;\|x^{\ast}\|\ls 1\}$.
We may assume that $\E \|\sum_{i\gs 1} y_i\va_i\|<\infty$ since otherwise there is nothing to prove.
Consequently for each $x^{\ast}\in B^{\ast}$ series $\sum_{i\gs 1}x^{\ast}(y_i)\va_i$ is convergent
which is equivalent to $\sum_{i\gs 1}(x^{\ast}(y_i))^2<\infty$.
Let $Q:B^{\ast}\ra \ell^2$ be defined by
$Q(x^{\ast})=(x^{\ast}(y_i))_{i\gs 1}$. It is clear that $Q:B^{\ast}/\ker Q\ra \ell^2$ is a linear
isomorphism on the closed linear subspace of $\ell^2$. We apply Theorem \ref{thm:2} to get the following
result.
\begin{coro}
Suppose that $Q$ is onto $\ell^2$ then (\ref{Ole1}) implies (\ref{Ole0}). 
\end{coro}
Unfortunately if $Q$ is not onto $\ell^2$ then the above argument fails. Still it is believed that
the comparison holds. A partial result can be deduced from Theorem \ref{thm:1} namely
\begin{coro}
Suppose that for each $x^{\ast}\in B^{\ast}$ and $p\gs 0$
\be\label{ole6}
\inf_{|I^c|\ls Cp}\sum_{i\in I}|x^{\ast}(x_i)|^2\ls C^2\inf_{|I^c|\ls p}\sum_{i\in I} |x^{\ast}(y_i)|^2.
\ee
Then (\ref{Ole0}) holds, i.e. 
$$
\E\|\sum_{i\gs 1}x_i\va_i\|\ls K\E\|\sum_{i\gs 1}y_i\va_i\|.
$$
\end{coro} 
\begin{dwd}
It suffices to notice that (\ref{ole6}) implies (\ref{ole2}) and then apply Theorem \ref{thm:1}.
\end{dwd}

\section*{Acknowledgments}
We would like to thank prof. Kwapie{\'n} for comments on the shape of this paper and helpful discussion about Theorem \ref{theo2}.

\end{document}